\documentclass[12pt, oneside]{article}
\usepackage{fontspec}
\usepackage{geometry, comment}
\usepackage[parfill]{parskip}
\usepackage{graphicx}
\usepackage[normalem]{ulem}

\usepackage{amsmath, amssymb, amsthm, mathrsfs, bbm}
\usepackage[hidelinks]{hyperref}
\usepackage{mathtools, xcolor, centernot, bm, setspace}

\usepackage[bottom]{footmisc}
\usepackage{enumerate}
\allowdisplaybreaks

\usepackage{titlesec}
\titleformat{\section}
{\centering\normalfont\bfseries}
{\thesection.}
{0.5em}
{\MakeUppercase}

\titleformat{\subsection}
{\normalfont\bfseries}
{\thesubsection.}
{0.5em}
{}

\usepackage{abstract}
\newcommand*\norm[1]{
	\vert\vert{#1}\vert\vert
}

\usepackage{tikzlings}
\usepackage{xparse}
\usepackage{etoolbox}

\newcommand{\A}{\mathcal{A}}
\newcommand{\B}{\mathcal{B}}
\newcommand{\E}{\mathbb{E}}
\newcommand{\F}{\mathfrak{F}}
\newcommand{\G}{\mathfrak{S}}

\newcommand{\bbN}{\mathbb{N}}
\renewcommand{\P}{\mathbb{P}}
\newcommand{\R}{\mathbb{R}}
\newcommand{\emptyword}{\varnothing}

\newcommand{\abs}[1]{\vert #1 \vert}

\renewcommand{\epsilon}{\varepsilon}

\newtheorem{thm}{Theorem}[section]
\newtheorem{cor}[thm]{Corollary}
\newtheorem{lemma}[thm]{Lemma}

\newtheorem{prop}[thm]{Proposition}

\theoremstyle{definition} 
\newtheorem{defn}[thm]{Definition}
\newtheorem{remark}{Remark}

\AtEndEnvironment{defn}{\hfill\ensuremath{\blacksquare}}
\AtEndEnvironment{ex}{\hfill\ensuremath{\blacksquare}}

\title{\bf \Large On the dimension of $V$-variable fractals}
\author{\normalsize Corey de Wit}
\date{}

\begin{document}
	\maketitle
	
	\begin{abstract}
		Random recursive sets have been used to model a wide range of fractal phenomena, but they can be slow and difficult to generate. $V$-variability address this problem by allowing these sets at most $V$ distinct shapes at each level of magnification. The formulae for the almost sure Hausdorff dimension of both classes is known, and it has been conjectured that (under mild assumptions) the $V$-variable dimension converges to the random recursive dimension as $V\to\infty$. We prove this conjecture. \let\thefootnote\relax\footnote{\emph{Date}: \today \\ \indent\hspace{8pt} The author would like to thank Jonathan Fraser for conversations regarding the problem, Antti K{\" a}enm{\" a}ki for an inspiring minicourse on fractal geometric applications of subadditivity, and Johannes Livengood van Vliet for revisions. The author was financially supported by a PhD studentship provided by the University of St Andrews in associated with the EPSRC Open Fellowship (EP/Z533440/1) held by Jonathan Fraser.}
	\end{abstract}
	
	\section{Introduction}
	
	Iterated function systems (IFSs) are a popular method for constructing fractals exhibiting self-similarity, and there is a rich and growing literature studying their attractor's dynamics and dimension. Random iterated function systems generalise classical IFSs by introducing randomness into the construction of the attractor over a family of IFSs. Well studied cases include the random homogeneous model, whereby a single IFS is randomly chosen at each level and applied uniformly across that level, and the random recursive model, where an IFS is chosen independently for each piece at every level of the construction (see \cite{falconer, mauldin-williams}). $V$-variable sets were introduced by Barnsley, Hutchinson, and Stenflo to interpolate between these two models by requiring there to be at most $V$ distinct local patterns at each level of magnification \cite{premain, main}, with $V=1$ the random homogeneous case and $V=\infty$ the random recursive case.
	
	For random systems made up of similarities and satisfying an open set condition, the Hausdorff dimension of almost every random homogeneous and recursive set can be readily calculated as a probabilistic version of the classical similarity dimension. The case for $V$-variable sets however relies on calculating the Lyapunov exponent of the decision trees \cite{main-dim}. Methods for such dimension calculations have since been researched and extended to more general settings, see for example \cite{code-trees,sascha2}. Nevertheless, it has remained a folklore conjecture since 2012, supported by the computational results in \cite{main-dim}, that the dimension of a typical $V$-variable set approaches the dimension of a typical random recursive set as $V\to\infty$. This is explicitly formulated in {\cite{sascha}}, which further conjectures the almost sure dimension as a function of $V$ is strictly increasing (provided its endpoint values are distinct). We prove the first of these conjectures by piecing together techniques in \cite{main,main-dim} and \cite{mauldin-williams-dim}.
	
	\pagebreak
	\section{$V$-variable fractals}
	
	We will first precisely define the construction of $V$-variable iterated function systems.
	
	\subsection{Random Iterated Function Systems \label{subsec:rifs}}
	
	\begin{defn}\label{def:rifs}
		Let $F = \{F_\lambda = \{f_{i,\lambda} \} \}_{\lambda\in\Lambda}$ be a family of IFSs consisting of contracting similitudes on $\R^d$, each with $N_\lambda$ maps and associated contraction factors $r_{i,\lambda}$. Let $\mu$ be a probability measure on the (possibly infinite) indexing set $\Lambda$. We say the pair $(F,\mu)$ is a \emph{random iterated function system} (RIFS). We make the following assumptions:
		\[ \inf N_\lambda \geq 2, \quad N := \sup N_\lambda < \infty, \]\[ 0 < r_{\min} := \inf r_{i,\lambda} \leq \sup r_{i,\lambda} =: r_{\max} < 1, \quad \sup \abs{f_{i,\lambda}(0)-0} < \infty . \]
	\end{defn}
	
	As mentioned in the introduction, the goal of an RIFS is to apply different function sets at each iteration depth. The choices of such IFSs can be written as a tree. Let $[N]^*$ be the set of finite words (including the empty word $\emptyword$) where $[n] := \{1,2,\ldots,n\}$, let $[N]^n$ denote the length $n$ words, and let $\abs{v}$ denote the length of $v\in[N]^*$. For $0\leq m<n$ and $v\in[N]^n$, we write $v\vert m = v_1v_2\cdots v_m$ for truncation, and for $v,w\in[N]^*$ we will write $vw$ for their concatenation. We will also view $[N]^*$ as the rooted $N$-ary tree, where $\emptyword$ is the root, the children of any $v\in[N]^*$ are $v1, v2, \ldots, vN$, and $[N]^n$ is the $n$-th level of the tree. A subtree of $[N]^*$ is a nonempty subset $\tau\subset [N]^*$ such that $v\in \tau$ implies $v \vert k \in \tau$ for all $k<\abs{v}$. We will write $\tau_n$ for the $n$-th level of $\tau$ and $\tau\vert_v$ for the subtree rooted at $v$. We will consider coloured subtrees of $[N]^*$ defined inductively by associating $\Lambda$-colourings to nodes at each level.
	
	\begin{defn}\label{def:tree-orig}
		We say $\tau$ is a $\Lambda$-coloured tree for an RIFS $(F,\mu)$ if it is a subtree of $[N]^*$ constructed in the following way: for every $n\in\bbN_0$ and node $v$ at level $n$ of the partially constructed tree,
		\begin{enumerate}[(1)]
			\item Colour $v$ with some $\lambda\in\Lambda$, 
			\item Build level $n+1$ nodes $v1,v2,\ldots, vN_\lambda$.
		\end{enumerate}
		Abusing notation, let $\tau$ denote both the constructed coloured tree as well as the colouring map which determined it (i.e. the partial function $\tau : [N]^* \rightharpoonup \Lambda$). Let $\Omega$ denote the collection of such trees. For $k\in\bbN$ and $v,w\in\tau_k$, we say the two coloured trees $\tau\vert_v$ and $\tau\vert_w$ are isomorphic if $vu\in\tau \iff wu\in\tau$ and $\tau(vu) = \tau(wu)$ for every $u\in [N]^*$, and say $\tau$ is \emph{$V$-variable} if there are at most $V$ isomorphism classes of coloured subtrees rooted at level $k$, for every $k\in\bbN_0$. Let $\Omega_V$ denote the set of such trees, so $\Omega_\infty = \Omega$. 
	\end{defn}
	
	\begin{defn}
		The attractor for the RIFS $(F,\mu)$ with respect to $\tau$ is 
		\begin{equation}\label{eq:attractor}
			F_\tau := \lim_{n\to\infty} \bigcup_{v\in\tau_n} f_{v_1, \tau(\emptyset)}\circ f_{v_2, \tau(v_1)}\circ\cdots\circ f_{v_n, \tau(v\vert n-1)}(K) 
		\end{equation}
		with respect to the Hausdorff metric, independent of the choice of non-empty compact set $K\subseteq\R^d$. We assume $F$ satisfies the \emph{Uniform Open Set Condition} (UOSC), that is there exists an open set $U\subseteq \R^d$ such that $f_{i,\lambda}(U)\subseteq U$ and $f_{i,\lambda}(U)\cap f_{j,\lambda}(U) = \emptyset$ for every $\lambda\in\Lambda$ and $1\leq i\neq j \leq N_\lambda$. Further, set 
		\[ \G_\lambda^s = \sum_{i=1}^{N_\lambda} r_{i,\lambda}^s , \quad  \G^s := \E_{\lambda\sim\mu} (\G_\lambda^s) , \quad \underline{\G}^s := \inf_\lambda\{\G_\lambda^s\}, \quad \overline{\G}^s := \sup_\lambda\{\G_\lambda^s\} \]
		Note our assumption on $N_\lambda \geq 2$ implies $\E_{\lambda\sim\mu} (\log \G^0_\lambda) > 0$ and $\G^0 > 1$.
	\end{defn}
	
	\begin{thm}\cite{falconer, mauldin-williams-dim}
		The Hausdorff dimension of a typical random homogenous attractor is the unique $s$ such that $\E_{\lambda\sim\mu}(\log \G^s_\lambda) = 0$, and the Hausdorff dimension of a typical random recursive attractor is the unique $s$ such that $\G^s = 1$. 
	\end{thm}
	
	\subsection{Another perspective of trees}
	
	Definition \ref{def:tree-orig} constructs arbitrarily coloured trees, but does not impose $V$-variability. Our aim in this section is to introduce an alternative tree construction which guarantees the property by tracking the isomorphism class information of each subtree at each level and encoding it in their root vertices. To motivate a formal definition, we will demonstrate a conversion procedure using a tree constructed by our current definition.
	
	Choose some $\tau\in\Omega_V$. Fix some $\lambda_0\in\Lambda$ and $i_0\in[V]$, these will act as `default' choices for colour and isomorphism class. For every $n\in\bbN_0$, let $V_n$ denote the number of isomorphism classes of subtrees at level $n$, index them $1,2,\ldots, V_n$, and let $v^{(n)}_i$ denote the root word of the $i$-th isomorphism class representative. Now build the $V$-tuple  
	\[ \lambda^{(n)} = \Big(\tau(v^{(n)}_1), \tau(v^{(n)}_2), \ldots, \tau(v^{(n)}_{V_n}), \underbrace{\lambda_0,\ldots,\lambda_0}_{V-V_n} \Big) , \]
	and a $V\times N$ matrix $J^{(n)}$, setting its $ij$-th element to be the unique $k$ such that $\tau\vert_{v^{(n)}_i j} = \tau\vert_{v^{(n+1)}_k}$ for $j \leq N_{\lambda_i^{(n)}}$ and $i_0$ otherwise, i.e. the index of the isomorphism class which the $j$-th child of $v^{(n)}_i$ belongs to. So the tuple $a_n = (\lambda^{(n)}, J^{(n)})$ records the root node colour of each subtree at the $n$-th level, and determines the isomorphism class index of each subtree at the $(n+1)$-th level in $\tau$'s construction. Moreover, the sequence $(a_n)_{n\in\bbN}$ determines the tree. One should note we assigned constants in redundant situations, and chose a specific ordering of the isomorphism classes at each level, so the assignment of the sequence $(a_n)_{n\in\bbN_0}$ to $\tau$ is many-to-one. We remark this information can be turned into a $V$-tuple of attractors (see \cite{superfractals}), but that isn't relevant to the work in this paper. 
	
	\begin{defn}\label{def:a-tree}
		Let $\A_V = \Lambda^V\times[V]^{V\times N}$, and fix some $i\in[V]$ and 
		\[ a = (\lambda^{(0)},J^{(0)})(\lambda^{(1)},J^{(1)})\cdots\in\A_V^{\bbN_0}, \quad\text{where }\lambda^{(n)} = (\lambda^{(n)}_1,\lambda^{(n)}_2,\cdots, \lambda^{(n)}_V), J^{(n)} = (J^{(n)})_{ij}. \] 
		Define the $\Lambda$-coloured tree $\tau^{a,i}$ and isomorphism class indexing function $\iota:\tau^{a,i}\to[V]$ the following way: for each $n\in\bbN_0$, each $v\in [N]^n$ will already be assigned some index $\iota(v)$ (with initialisation $\iota(\emptyword)=i$), so
		\begin{enumerate}[(1)]
			\item For each $v$ at level $n$ of the partially constructed tree, 
			\begin{enumerate}[(i)]
				\item Set the colour $\tau^{a,i}(v)$ of $v$ to $\lambda_{\iota(v)}^{(n)}$ and 
				\item Build nodes $v1,v2,\ldots, vN_{\tau^{a,i}(v)}$
			\end{enumerate}
			\item For each $v \in [N]^n$ and $j\in[N]$, let $\iota(vj) = (J^{(n)})_{\iota(v)j}$. 
		\end{enumerate}
		We will write $\tau^a := \tau^{a,1}$ (i.e. giving the root node the sensible isomorphism class indexing). With $U_V$ the uniform probability measure on $[V]$ for $V\in\bbN$, let
		\[ \P_V^\mu := (\mu^{\otimes V} \otimes U_V^{\otimes VN})^{\otimes\bbN_0} \]
		be the canonical measure (with respect to $\mu$ and $U_V$) on the product space $\A_V^{\bbN_0}$. Note $\P_V^\mu$ is invariant under $\sigma$ and permutation of $[V]$. Finally, let $\B_V = \B(\A_V^{\bbN_0})$ be the Borel $\sigma$-algebra generated by cylinder sets, and $\B_{n,V}\subset \B_V$ be the $\sigma$-algebra generated by cylinder sets of words up to length $n$ (i.e. the smallest $\sigma$-algebra which makes projections onto the $k$-th coordinate of the string $a_0a_1\ldots$, for $k<n$, measurable).
	\end{defn}
	\begin{remark}
		We will similarly let the probability measure on $\Omega$, denoted $\P_\infty^\mu$, be the pushforward of $\mu^{\otimes [N]^*}$ under the mapping $\Lambda^{[N]^*}\to\Omega$. Let $\B_\infty$ be the $\sigma$-algebra generated by sets $\{\tau\in\Omega \mid v\in\tau, \tau(v)\in \Lambda' \}$ for $v\in[N]^*$ and $\mu$-measurable $\Lambda'\subset\Lambda$, and $\B_{n,\infty}\subset\B_\infty$ be the $\sigma$-algebra generated by such sets with $\abs{v}<n$. Additionally one can note that, denoting $\sigma$ the left shift on $\A_V^{\bbN_0}$, for $v\in\tau_n^{a,i}$
		\begin{equation}\label{eq:restriction-is-shift}
			\tau^{a,i}\vert_v = \tau^{\sigma^n(a),\iota(v)} .
		\end{equation}
	\end{remark}
		
	\subsection{Flows on trees}
	We can now define the flow of IFS $s$-mass between two subtrees rooted at different levels belonging to different isomorphism classes and use it to calculate the dimension of $V$-variable sets.
	\begin{defn}
		The $V\times V$ flow matrix $M_a$ for some $a = (\lambda, J)\in\A_V$ is given by
		\begin{equation}\label{eq:flowmatrix}
			(M_a(s))_{ij} = \sum_{k=1}^{N_{\lambda_i}} \{r_{k,\lambda_i}^s \mid J_{ik} = j \} , 
		\end{equation}
		and for $a\in\A_V^{\bbN_0}$ write $M_{a\vert n}(s) = M_{a_0}(s)M_{a_1}(s) \cdots M_{a_{n-1}}(s)$. Let $\norm{M} = \sum_{i,j} M_{i,j}$. 
	\end{defn}
	\begin{thm}[{\cite[Main Theorem]{main-dim}}]\label{thm:main-dim}
		Let $(F,\mu)$ be a RIFS. The Lyapunov exponent of a $V$-variable system is equal to 
		\begin{equation}\label{eq:lyapunov}
			\gamma_V(s) = \lim_{n\to\infty} \E_{a\sim\P^\mu_V} \Bigg(\frac{1}{n} \log \norm{M_{a\vert n}(s)}\Bigg) ,
		\end{equation}
		and is monotonically decreasing in $s$. Additionally, for almost every $a\in\A_V^{\bbN_0}$, the Hausdorff dimension of $F_{\tau^a}$ is the unique zero of $\gamma_V(s)$, so define the function $D: \bbN \cup\{\infty\} \to \R^+$ by
		\[ D(V) = \mbox{a.s.-}\dim_H(F_{\tau^a}) , \quad D(\infty) = \mbox{the unique $s$ such that $\G^s=1$.}  \]
	\end{thm}
	
	\pagebreak
	\section{Proof of the dimension limit}
	
	\subsection{Properties of flows}
	
	In \cite[\S 5]{main-dim}, they look at the sum over the $i$-th row of the finite product of matrices in the definition of $\gamma_V(s)$, however this direct definition can only be discussed in the $V$-variable setting. We will look at an object studied in \cite{mauldin-williams-dim} for the $V=\infty$ case which extends the aforementioned notion and contains results similar to Lemma \ref{lem:flow-recursive} and Proposition \ref{prop:flow-Ktoinfty-bound}. For $K\in\bbN_0$, $s\geq0$ and $\tau\in\Omega$, let
	\begin{equation}\label{eq:flow}
		\F(K,s,\tau) := \sum_{v\in\tau_K} \prod_{k<K} r_{v_{k+1},\tau(v\vert k)}^s 
	\end{equation}
	be the total flow of IFS $s$-mass over the first $K$ levels of $\tau$. We will first note the recursive nature of this definition, then relate this flow to the product of matrices in Equation \eqref{eq:lyapunov}.
	
	\begin{lemma}\label{lem:flow-recursive}
		For all $K\in\bbN_0$, $s\geq 0$ and $\tau\in\Omega$,
		\begin{equation}\label{eq:flow-recursive-1}
			\F(K+1,s,\tau) = \sum_{i=1}^{N_{\tau(\emptyword)} }r_{i,\tau(\emptyword)}^s \F(K,s,\tau\vert_i)
		\end{equation}
		and
		\begin{equation}\label{eq:flow-recursive-2}
			(\underline{\G}^s)^K \leq \F(K,s,\tau)\leq (\overline{\G}^s)^K .
		\end{equation}
	\end{lemma}
	\begin{proof}
		We can isolate the first or last term in the product in Eq.~\eqref{eq:flow}. The former gives
		\begin{equation}\label{eq:break-off-first}
			\F(K+1,s,\tau) = \sum_{v\in\tau_{K+1}} r_{v_1,\tau(\emptyword)}^s \prod_{1\leq k<K+1} r_{v_{k+1},\tau(v\vert k)}^s , 
		\end{equation}
		which is equivalent to Eq.~\eqref{eq:flow-recursive-1}. The latter gives
		\[ \F(K,s,\tau) = \sum_{v\in\tau_{K}} r_{v_{K},\tau(v\vert K-1)}^s \prod_{k<K-1} r_{v_{k+1},\tau(v\vert k)}^s , \]
		thus 
		\[ \underline{\G}^s\F(K-1,s,\tau) \leq \F(K,s,\tau)\leq \overline{\G}^s\F(K-1,s,\tau) . \]
		Iterating the above gives Eq.~\eqref{eq:flow-recursive-2}.
	\end{proof}
	
	\begin{prop}\label{prop:local-window-flow}
		For any $a\in\A_V^{\bbN_0}$, $K\in\bbN_0$, $s\geq 0$ and $i\in[V]$,
		\[ (M_{a\vert K}(s){\bf1})_i = \F(K,s,\tau^{a,i}) . \]
	\end{prop}
	\begin{proof}
		We will prove this inductively over $K$. With the base case trivial, assume true for a fixed $K$, and denote $a_0 = (\lambda, J)$. Then
		\begin{align*}
			(M_{a_0}(s)M_{\sigma(a)\vert K}(s){\bf1})_i &= \sum_j (M_{a_0}(s))_{ij} \F(K,s,\tau^{\sigma(a),j}) \\
			&= \sum_j \F(K,s,\tau^{\sigma(a),j}) \sum_{k=1}^{N_{\lambda_i}} \{ r_{k,\lambda_i}^s \mid J_{ik} = j \} \\
			&= \sum_{k=1}^{N_{\lambda_i}} r_{k,\lambda_i}^s \F(K,s,\tau^{\sigma(a),J_{ik}}) . 
		\end{align*}
		Since $\tau^{a,i}(\emptyword) = \lambda_i$ by the first step of tree construction in Definition \ref{def:a-tree}, and $\tau^{a,i}\vert_k = \tau^{\sigma(a), J_{ik}}$ by Eq.~\eqref{eq:restriction-is-shift}, then Eq.~\eqref{eq:flow-recursive-1} gives us the desired equality.
	\end{proof}
	
	This has two immediate consequences. 
	\begin{cor}\label{cor:lyapunov-is-flow}
		For any $a\in\A_V^{\bbN_0}$, $s\geq 0$ and $K \geq 1$,
		\[ \norm{M_{a\vert K}(s)} =  \sum_{i=1}^V \F(K,s,\tau^{a,i}) . \]
	\end{cor}
	
	\begin{cor}\label{cor:expected-flow}
		For every $s\geq 0$, $K\geq 1$ and $1\leq i\leq V$, 
		\[ \E_{a\sim \P^\mu_V} \Big(\F(K,s,\tau^{a,i})\Big) = (\G^s)^K = \E_{\tau\sim \P^\mu_\infty} \Big(\F(K,s,\tau)\Big) \]
	\end{cor}
	\begin{proof}
		Since each $a_n$ are independent, then 
		\[ \E_{a\sim\P^\mu_V} \Big(M_{a\vert K}(s)\Big)_{ij} = \bigg(\E_{a\sim \mu^{\otimes V}\otimes U_V^{\otimes VN}} \Big(M_a(s)\Big)^K \bigg)_{ij}.\]
		Additionally, $(M_a(s){\bf1})_i = \G_{\lambda_i}^s$ for any $a\in\A_V$, so Proposition \ref{prop:local-window-flow} gives the first equality. Additionally, since the colouring of the root is independent to the colouring of the rest of the tree with respect to $\P_\infty^\mu$, we get the second equality by taking the expectation of both sides in Eq.~\eqref{eq:flow-recursive-1} and iterating.
	\end{proof}
	
	\subsection{Expectation results}
	
	By Corollary \ref{cor:lyapunov-is-flow}, it suffices to work with these flows when studying $\gamma_V(s)$. We will first adapt \cite[Theorem 49]{main} to our setting of flows. For a fixed $s$ we will write
	\[ \Phi_K(\tau) := \frac{1}{K}\log\F(K,s,\tau), \]
	and will write $\Phi(\tau)$ when $K$ is fixed. 
	
	\begin{prop}\label{prop:flow-Vtoinfty-bound}
		For all $K,V\in\bbN$ and $s\geq 0$
		\begin{equation}\label{eq:flow-Vtoinfty-bound}
			\bigg\vert \E_{a\sim\P_V^\mu} \Big(\Phi_K(\tau^a)\Big) - \E_{\tau\sim\P_\infty^\mu} \Big(\Phi_K(\tau)\Big) \bigg\vert \leq \frac{N^{2K} }{V}\log\frac{\overline{\G}^s}{\underline{\G}^s} 
		\end{equation}
	\end{prop}
	\begin{proof}
		Fix $K$, $V$ and $s$. For any $a\in\A_V^{\bbN_0}$ and $\tau\in\Omega$, by Eq.~\eqref{eq:flow-recursive-2}
		\[ \bigg(\frac{\overline{\G}^s}{\underline{\G}^s}\bigg)^{-K} \leq \frac{\F(K,s,\tau^a)}{\F(K,s,\tau)} \leq \bigg(\frac{\overline{\G}^s}{\underline{\G}^s}\bigg)^K , \]
		so taking the logarithm of both sides and dividing by $K$ gives the general bound
		\begin{equation}\label{eq:flow-Vtoinfty-prelim-bound}
			\vert \Phi(\tau^a) - \Phi(\tau) \vert \leq \log\frac{\overline{\G}^s}{\underline{\G}^s} .
		\end{equation}
		For $N^K\geq V$ this immediately implies Eq.~\eqref{eq:flow-Vtoinfty-bound}, so assume $N^K<V$. Then there exist $a\in\A_V^{\bbN_0}$ such that $\iota$ is injective on $[N]^k$ for $k<K$; let $E\in\B_{K,V}$ denote the collection of such $a$. For all $k\leq K$, note: (a) $\iota\vert_{[N]^k}$ is only determined by $J^{(l)}$ for $l<k$ and so is independent of $\lambda$, and (b) the injectivity of $\iota\vert_{[N]^k}$ implies the colouring $\lambda^{(k)}_{\iota(v)}$ of each $v\in [N]^k$ is chosen i.i.d.. Since $\F(K,s,\tau)$ only sees the first $K$ levels of $\tau$, then
		\begin{equation}\label{eq:injective-event}
			\E_{a\sim\P_V^\mu} \Big(\Phi(\tau^a) \mid E \Big) = \E_{\tau\sim\P_\infty^\mu} \Big(\Phi(\tau)\Big) , 
		\end{equation}
		So expanding the left hand side of Eq.~\eqref{eq:flow-Vtoinfty-bound} and applying Eq.~\eqref{eq:injective-event}, we get
		\begin{align*}
			\bigg\vert \E_{a\sim\P_V^\mu} \Big(\Phi(\tau^a) \mathbbm{1}_{E^c} &\Big) + \bigg[\E_{a\sim\P_V^\mu} \Big(\Phi(\tau^a) \mathbbm{1}_E \Big) - \P_V^\mu(E)\cdot\E_{\tau\sim\P_\infty^\mu} \Big(\Phi(\tau)\Big) \bigg] - \P_V^\mu(E^c)\cdot\E_{\tau\sim\P_\infty^\mu} \Big(\Phi(\tau)\Big) \bigg\vert \\
			&= \bigg\vert \E_{a\sim\P_V^\mu} \bigg(\bigg[\Phi(\tau^a) - \E_{\tau\sim\P_\infty^\mu} \Big(\Phi(\tau)\Big) \bigg] \mathbbm{1}_{E^c} \bigg) \bigg\vert \\
			&\leq \P_V^\mu(E^c) \sup_{a\in\A_V^{\bbN_0}} \bigg\vert \Phi(\tau^a) - \E_{\tau\sim\P_\infty^\mu} \Big(\Phi(\tau)\Big) \bigg\vert \\
			&\leq \P_V^\mu(E^c) \log\frac{\overline{\G}^s}{\underline{\G}^s} ,
		\end{align*}
		where the second inequality is a result of applying the general bound Eq.~\eqref{eq:flow-Vtoinfty-prelim-bound}. Finally, for an $a\in E^c$ there is a least $k<K$ such that two vertices at level $k$ are labelled the same isomorphism class index. By minimality of $k$, this is equivalent to there being two distinct entries of $J^{(k-1)}$ which have equal value. Since entries are chosen i.i.d., the probability a fixed pair agrees is $1/V$, and since there are fewer than ${N^k\choose2}$ choices of pairs at each level, then 
		\[ \P_V^\mu(E^c) \leq \frac{1}{V}\sum_{k=0}^{K-1} {N^k\choose2} \leq \frac{N^{2K}}{V} . \]
	\end{proof}
	
	Thus the difference between the expected log-normed flow to a fixed level of a $V$-variable and random tree shrinks to zero as $V\to\infty$. We now want to find a Lyapunov exponent representation for RIFS systems equivalent to Eq.~\eqref{eq:lyapunov}.
	
	\begin{prop}\label{prop:flow-Ktoinfty-bound}
		For any $s$ such that $\G^{2s}<(\G^s)^2$,
		\[ \E_{\tau\sim\P_\infty^\mu} \Big(\Phi_K(\tau)\Big) \xrightarrow{K\to\infty} \log(\G^s) . \]
	\end{prop}
	\begin{proof}
		Let $\widetilde{\F}(K,s,\tau) := \F(K,s,\tau)/(\G^s)^K$ denote the normalised flow. We claim its second moment is uniformly bounded in $K$. By \eqref{eq:flow-recursive-1} and the independence of the colouring of the root to the rest of the tree with respect to $\P_\infty^\mu$, 
		\begin{align*}
			\E_{\tau\sim\P_\infty^\mu}\Big(\widetilde{\F}(K+1,s,\tau)^2\Big) &= \frac{1}{(\G^s)^{2K+2}}\E_{\tau\sim\P_\infty^\mu}\Bigg(\bigg( \sum_{i=1}^{N_{\tau(\emptyword)}} r_{i,\tau(\emptyword)}^s \F(K,s,\tau\vert_i) \bigg)^2 \Bigg) \\
			&= \frac{1}{(\G^s)^{2K+2}}\sum_{i=1}^N \E_{\lambda\sim\mu} \Big( r_{i,\lambda}^{2s} \Big) \E_{\tau\sim\P_\infty^\mu}\Big( \F(K,s,\tau)^2 \Big)  \\
			&\qquad+ \frac{1}{(\G^s)^{2K+2}} \sum_{i\neq j} \E_{\lambda\sim\mu} \Big( r_{i,\lambda}^{s} r_{j,\lambda}^{s} \Big) \E_{\tau\sim\P_\infty^\mu} \Big( \F(K,s,\tau) \Big)^2  \\
			&= \frac{1}{(\G^s)^2} \G^{2s} \E_{\tau\sim\P_\infty^\mu}\Big( \widetilde{\F}(K,s,\tau)^2 \Big) \\
			&\qquad+ \frac{1}{(\G^s)^2} \Bigg(  \sum_{i,j} \E_{\lambda\sim\mu} \Big( r_{i,\lambda}^{s} r_{j,\lambda}^{s} \Big) - \sum_{i=1}^N \E_{\lambda\sim\mu} \Big(r_{i,\lambda}^{2s} \Big)  \Bigg)   \\
			&= \frac{\G^{2s}}{(\G^s)^2} \E_{\tau\sim\P_\infty^\mu}\Big( \widetilde{\F}(K,s,\tau)^2 \Big) + \frac{\E_{\lambda\sim\mu}\Big((\G_\lambda^s)^2\Big) - \G^{2s}}{(\G^s)^2} ,
		\end{align*}
		where we let $r_{i,\lambda} = 0$ for $N_\lambda<i\leq N$, and use Corollary \ref{cor:expected-flow} for the third equality. Also note $\widetilde{\F}(1,s,\tau) < \infty$ for every $\tau\in\Omega$, so by our assumption on $s$ this completes the claim, and moreover shows $\sup_K \norm{\widetilde{\F}(K,s,\tau)}_2 < \infty$. Now notice that $\E_{\tau\sim \P^\mu_\infty}(\widetilde{\F}(K,s,\tau)) = 1$ for every $K$, and further
		\[ \E_{\tau\sim \P^\mu_\infty}\Big(\widetilde{\F}(K+1,s,\tau) \mid \B_{K,\infty} \Big) = \widetilde{\F}(K,s,\tau), \]
		so the sequence $(\tau\mapsto\widetilde{\F}(K,s,\tau))_{K\in\bbN}$ is a martingale adapted to $(\B_{K,\infty})_{K\in\bbN}$, and converges $\P_\infty^\mu$-a.e. and in $L^2$, to some $\widetilde{\F}(s,\tau)$ by the Martingale Convergence Theorem \cite{martingale}. 
		
		This has two consequences. First, the convergence paired with uniform boundedness implies that $\E_{\tau\sim \P^\mu_\infty}(\widetilde{\F}(s,\tau)) = 1$. Second, dividing both sides of Eq.~\eqref{eq:flow-recursive-1} by $(\G^s)^{K+1}$ and taking the limit as $K\to\infty$ gives 
		\[ \widetilde{\F}(s,\tau) = \frac{1}{\G^s} \sum_{i=1}^{N_{\tau(\emptyword)}}r_{i,\tau(\emptyword)}^s \widetilde{\F}(s,\tau\vert_i) \quad \P_\infty^\mu\text{-a.e.,} \]
		so $\widetilde{\F}(s,\tau) = 0$ if and only if $\widetilde{\F}(s,\tau\vert_i) = 0$ for all $1\leq i\leq N_{\tau(\emptyword)}$. Moreover,
		\[ \P_\infty^\mu(\widetilde{\F}(s,\tau)=0) = \E_{\lambda\sim\mu}\Big(\P_\infty^\mu(\widetilde{\F}(s,\tau)=0)^{N_\lambda} \Big) . \]
		By our assumption that $N_\lambda\geq2$, then $\P_\infty^\mu(\widetilde{\F}(s,\tau)=0) \in\{0,1\}$. The first consequence rules out the latter option, so $\widetilde{\F}(s,\tau)>0$ almost surely, and $\log\widetilde{\F}(s,\tau)\in(-\infty,\infty)$. Thus 
		\[ \Phi_K(\tau) = \log\G^s + \frac{1}{K}\log\widetilde{\F}(K,s,\tau) \xrightarrow{K\to\infty} \log\G^s \quad \P_\infty^\mu\text{-a.e.}. \]
		Finally by Eq.~\eqref{eq:flow-recursive-2}, $\abs{\Phi_K(\tau)}$ is bounded by a constant, so by the Dominated Convergence Theorem
		\[ \lim_{K\to\infty}\E_{\tau\sim\P_\infty^\mu}\Big(\Phi_K(\tau)\Big) = \E_{\tau\sim\P_\infty^\mu}\Big(\lim_{K\to\infty}\Phi_K(\tau)\Big) = \log\G^s . \]
	\end{proof}
	
	\begin{cor}\label{cor:flow-toinfty}
		For any $s$ such that $\G^{2s}<(\G^s)^2$, setting $K(V) := \lfloor \frac{1}{4}\log_N V \rfloor$,
		\[ \E_{a\sim\P_V^\mu} \Big(\Phi_{K(V)}(\tau^a)\Big) \xrightarrow{V\to\infty} \log(\G^s) . \]
	\end{cor}
	\begin{proof}
		Combine Propositions \ref{prop:flow-Vtoinfty-bound} and \ref{prop:flow-Ktoinfty-bound}, noting that $N^{2K(V)}/V \leq V^{-1/2}$.
	\end{proof}
	
	\subsection{Proof of $D(V) \to D(\infty)$}
	
	We will now bound $\gamma_V(s)$ by sequence and limit in Corollary \ref{cor:flow-toinfty}, allowing us to use a Squeeze theorem argument to prove the dimension limit. 
	
	\begin{prop}\label{prop:bounds}
		For any $K,V\in\bbN$ and $s\geq0$,
		\[ \E_{a\sim\P_V^\mu} \Big(\Phi_K(\tau^a)\Big) \leq \gamma_V(s) \leq \log \G^s . \]
	\end{prop}
	\begin{proof}
		By Jensen's inequality,
		\[ \E_{a\sim\P_V^\mu} \Big(\frac{1}{K}\log\sum_{i=1}^V \F(K,s,\tau^{a,i}) \Big) \leq \frac{1}{K}\log\E_{a\sim\P_V^\mu}\bigg(\sum_{i=1}^V \F(K,s,\tau^{a,i}) \bigg) = \frac{1}{K}\log V + \log\G^s, \]
		where the equality comes from Corollary \ref{cor:expected-flow}. Thus as $K\to\infty$ the first term vanishes and Corollary \ref{cor:lyapunov-is-flow} gives us the upper bound. 
		
		To show the lower bound, letting $x_K := \E_{a\sim\P_V^\mu}\log\F(K,s,\tau^a)$ we claim $(x_K)_{K\in\bbN}$ is superadditive. By isolating the $K$ terms of the product in Eq.~\eqref{eq:flow}, for any $\tau\in\Omega$ we get 
		\[ \F(L+K,s,\tau) = \sum_{v\in\tau_K} \bigg(\prod_{k<K}r_{v_{k+1},\tau(v\vert k)}^s \bigg) \F(L,s,\tau\vert_v) . \]
		Now let $w_v := \prod_{k<K}r_{v_{k+1},\tau^a(v\vert k)}^s / \F(K,s,\tau^a)$, so $\sum_{v\in\tau_K} w_v = 1$. Then by Jensen's inequality and Eq.~\eqref{eq:restriction-is-shift},  
		\begin{equation}\label{eq:subadditive-before-expectation}
			\log \F(L+K,s,\tau^a) \geq \log\F(K,s,\tau^a) + \sum_{v\in(\tau^a)_K} \sum_{i=1}^V w_v \log\F(L,s,\tau^{\sigma^K(a),i}) \mathbbm{1}_{\iota(v)=i} ,
		\end{equation}
		where the second term was achieved by splitting the possible evaluations of $\iota(v)$. By the independence of data between the first $K$ levels and the rest, as well as the invariance of $\P_V^\mu$ under permutation of $[V]$, we have
		\[ \E_{a\sim\P_V^\mu}\Big(w_v \log\F(L,s,\tau^{\sigma^K(a),i}) \mathbbm{1}_{\iota(v)=i} \Big) = x_L \cdot\E_{a\sim\P_V^\mu}\Big(w_v \mathbbm{1}_{\iota(v)=i} \Big) . \]
		Therefore, taking expectations with respect to $\P_V^\mu$ of both sides of \eqref{eq:subadditive-before-expectation} gives the claim. So by Fekete's lemma and Corollary \ref{cor:lyapunov-is-flow},
		\[ \E_{a\sim\P_V^\mu}\Big(\Phi_K(\tau^a)\Big) \leq \sup_K \frac{x_K}{K} = \lim_{K\to\infty}\frac{x_K}{K} \leq  \gamma_V(s) . \]
	\end{proof}
	
	In order to transfer convergence of $\gamma_V$ to its zeros, we require the following lemma.
	\begin{lemma}\label{lem:uniform-rate}
		For all $V\in\bbN$ and $s,t\geq 0$,
		\[ \gamma_V(s+t) \leq \gamma_V(s) - t\log r_{\max}^{-1} . \]
	\end{lemma}
	\begin{proof}
		By observation of \eqref{eq:flowmatrix}, for any $a\in\A_V$ we have $(M_a(s+t))_{ij} \leq r_{\max}^t (M_a(s))_{ij}$, thus
		\[ \norm{M_{a\vert n}(s+t)} \leq r_{\max}^{t n}\norm{M_{a\vert n}(s)} . \]
		Taking the logarithms and expectations of both sides, dividing by $n$, and letting $n\to\infty$ gives the result.
	\end{proof}
	\begin{cor}\label{cor:decreasing-lyapunov}
		For any $V\in\bbN$, $\gamma_V(s)$ is decreasing in $s$. Moreover, $D(1)\leq D(V)\leq D(\infty)$.
	\end{cor}
	\begin{proof}
		For the second part, by Proposition \ref{prop:bounds}
		\[ \gamma_V(D(V)) \geq \E_{a\sim\P_V^\mu} \Big(\Phi_1(\tau^a)\Big) = \E_{\lambda\sim\mu} \Big(\log \G^{D(V)}_\lambda\Big) = \gamma_1(D(V)), \]
		and
		\[ \gamma_V(D(\infty)) \leq \log \G^{D(\infty)} = 0. \]
	\end{proof}
	
	\begin{thm}
		Let $(F,\mu)$ be a RIFS. Then $D(V) \xrightarrow{V\to\infty}D(\infty)$.
	\end{thm}
	\begin{proof}
		Let $s^* = D(\infty)$. By setting $s = D(V)$ and $t = s^*-D(V)$ (which is non-negative by Corollary \ref{cor:decreasing-lyapunov}), Lemma \ref{lem:uniform-rate} gives
		\[ s^* - D(V) \leq -\frac{\gamma_V(s^*)}{\log r_{\max}^{-1}} . \]
		Finally, by our assumption that $\G^0 > 1$ then $s^*>0$, so by definition of $s^*$ we have $\G^{2s^*}< \G^{s^*} = 1 = (\G^{s^*})^2$. Thus by Corollary \ref{cor:flow-toinfty} and Proposition \ref{prop:bounds}, the Squeeze theorem gives 
		\[ \gamma_V(s^*) \xrightarrow{V\to\infty} \log\G^{s^*} = 0 . \]
	\end{proof}
	
	\bibliography{v-variable-dim.bib}	
	
	\bigskip
	\noindent
	\textsc{C. de Wit, School of Mathematics and Statistics, University of St Andrews, St Andrews,
		KY16 9SS, United Kingdom.}
	
	\noindent
	\textit{Email address:} \texttt{cdw9@st-andrews.ac.uk}
\end{document}